\documentclass[11pt]{article} 
\usepackage{geometry}
\usepackage{amsmath,amssymb,amsfonts,amsthm} 
\usepackage[all]{xy}
\DeclareMathOperator{\Pf}{Pf}

\renewcommand{\AA}{\mathbb{A}}
\newcommand{\CC}{\mathbb{C}}

\newcommand{\FF}{\mathbb{F}}
\newcommand{\PP}{\mathbb{P}}
\newcommand{\QQ}{\mathbb{Q}}

\newcommand{\al}{\alpha}
\newcommand{\be}{\beta}

\newcommand{\si}{\sigma}

\newcommand{\Be}{\mathcal{B}}
\newcommand{\Ce}{\mathcal{C}}
\newcommand{\Oh}{\mathcal{O}}
\newcommand{\eX}{\mathcal{X}}

\newcommand{\Xtilde}{\widetilde X}

\newcommand{\Sbar}{\bar{S}}

\newcommand{\fie}{\varphi}
\newcommand{\two}{\textrm{II}}

\renewcommand{\qed}{\hfill\square}

\DeclareMathOperator{\Pic}{Pic}

\DeclareMathOperator{\Proj}{Proj}

\newtheorem{thm}{Theorem}[section]
\newtheorem*{thm-prime}{Theorem 1.3$^\prime$}
\newtheorem{conj}{Conjecture}[section]

\newtheorem{cor}{Corollary}[section]
\newtheorem{lemma}{Lemma}[section]
\newtheorem{rmk}{Remark}[section]
\newtheorem{dfn}{Definition}[section]
\newtheorem{exa}{Example}[section]
\newenvironment{pf}{\paragraph{Proof}}{\par\medskip}

\title{K3 transitions and canonical 3-folds}
\author{Stephen Coughlan}

\date{}

\begin{document}
\maketitle
\begin{abstract}
We introduce K3 transitions as a geometric approach to studying canonical 3-folds. These transitions link different deformation classes of canonical 3-folds via a combination of birational contractions and smoothings. As applications, we investigate some basic properties of the web of canonical 3-folds in small codimension and give an interesting example of a singularity in codimension 4 with obstructed smoothings, similar to the famous example of the affine cone over a degree 6 del Pezzo surface.
\end{abstract}
\section{Introduction}
A canonical 3-fold is a 3-dimensional variety $X$ with $K_X$ ample and at worst $\QQ$-factorial canonical singularities. Such varieties can be embedded in weighted projective space by taking the canonical model $\Proj R(X)$ of $X$, where $R(X)$ is the canonical graded ring
\[R(X,K_X)=\bigoplus_{n\ge0}H^0(X,nK_X).\]

In the 1980s, Iano-Fletcher \cite{Fletcher} and Reid \cite{YPG} produced lists of families of canonical 3-folds that are complete intersections in weighted projective space, and these lists were later proven to be complete by Chen, Chen and Chen \cite{Chen3}. Recently Brown, Kasprzyk and Zhou \cite{Brown-Kasprzyk-Zhu} discovered many new families that are not complete intersections, including 18 families in codimension 3 that are Pfaffian.

A family of canonical 3-folds is determined by the input data $\chi(\Oh_X)$, $P_2(X)$ and the basket of singularities $\mathcal{B}$, via the Riemann--Roch formula of \cite{YPG}, which computes the plurigenera $P_i(X)=h^0(nK_X)$ of $X$. The intricate combinatorics of such baskets has been heavily studied \cite{Chen-Chen}, and many restrictions on the possible baskets are known. Nevertheless, it remains difficult to determine whether a given triple $\chi,P_2,\mathcal{B}$ is realised as a canonical 3-fold.

Here we introduce \emph{K3 transitions}, transformations from one family of canonical 3-folds to another. The idea is to degenerate a canonical 3-fold $X$ to a 3-fold with a singularity that is locally the affine cone over a K3 surface. This is a particular example of a log canonical singularity. We resolve this singularity and take the canonical model $Y$ of the resolution. In good cases, $X$ and $Y$ are in different deformation families of canonical 3-folds.

K3 transitions can be used to traverse the geography of canonical 3-folds, keeping track of the changes that each transition makes to the data $\chi,P_2,\Be$. In this article, we study K3 transitions between the 137 families of canonical 3-folds of codimension $\le3$, studying what we call the \emph{web of canonical 3-folds}. 

The study of all log canonical singularities and their smoothings is very complicated. We know of other types of transitions that are relevant, including contractions of surfaces to a point or to a curve. We touch on some of these briefly, but mostly restrict our attention to cones over K3 surfaces.

Moreover, the moduli space of canonical 3-folds can have several irreducible components, and we could land in any one of them after performing a K3 transition. This does not happen in any of our examples in small codimension, but we return to this in \cite{Brown-Coughlan}.

We also exhibit an example of a singularity embedded in codimension 4 that has two topologically distinct smoothings, which will make an appearance in future work. Our obstructed singularity is distantly related to the cone over the degree 6 del Pezzo surface \cite{Altmann}.

Perhaps the strongest analogy for K3 transitions, is with primitive contractions on Calabi--Yau 3-folds \cite{Gross}, \cite{Kapustki}. There are other connections, such as projecting from points or lines on Fano 3-folds \cite{BKR} or primitive contractions on Calabi--Yau 3-folds \cite{Gross}, and elliptic singularities also play a prominent role in the study of geography of surfaces of general type \cite{Persson}.

\subsection{Elliptic singularities} The $3$-fold affine cone $\Ce_S$ over a K3 surface $S$ is called a \emph{simple elliptic 3-fold singularity} cf.~\cite[\S4]{C-3f}. For example, the cone over a quartic surface in $\PP^3$. An important point is that we do not require $S$ to be nonsingular; we allow affine cones over quasismooth K3 surfaces with $A_n$-singularities, embedded in weighted projective space $\PP(\al)=\PP(\al_0,\dots,\al_n)$. Quasismooth means precisely that the affine cone is nonsingular away from the vertex, i.e.~the singularity is isolated. An elliptic singularity $P$ on $X$ imposes one adjunction condition on $X$. That is, if $\Xtilde\to X$ is a resolution of $P$, then $|K_{\Xtilde}|$ is generated by those sections of $|K_X|$ which pass through $P$.

\subsection{K3 transitions}
Suppose we have a flat family $\eX$ whose general fibre $X$ is a canonical 3-fold, and whose central fibre $X_0$ has a (simple) elliptic singularity $P\in X_0$. Take the $\al$-weighted blow up $\si\colon Z\to X_0$ so that the exceptional divisor $S\subset Z$ is a quasismooth K3 surface $S\subset\PP(\al)$ with $A_n$-singularities.
\[\xymatrix
{&&Z\ar[dl]_{\si}\ar[dr]^{\pi}&\\
X\ar@{~>}[r]&X_0&&Y_0&Y\ar@{~>}[l]}\]
By the adjunction formula, $K_Z=\si^*K_{X}-S$. Thus $P$ is log canonical, and we can view $X\rightsquigarrow X_0$ as a degeneration to the boundary of the moduli space of canonical 3-folds. Define $Y_0=\Proj\bigoplus_{n=0}^\infty H^0(Z,nK_Z)$, the canonical model of $Z$. Then $Y_0$ is a canonical 3-fold containing the image $\Sbar$ of $S$ under $\pi$. In general, $\pi$ contracts several curves on $Z$ each to ordinary double points on $Y_0$ along $\Sbar$. Suppose there exists a partial smoothing of $Y_0$ to $Y$, a quasismooth canonical 3-fold. The procedure transforming $X$ to $Y$ is called a \emph{K3 transition in $S$} and we denote it by $X\to Y$.
\begin{exa}\label{ex!intro} The sextic 3-fold $X_6\subset\PP^4$ specialises to $X_0$ with a singularity $P$ which is locally the cone over a quartic surface. The projection $\pi_P\colon X\dashrightarrow\PP^3$ from $P$ to the complementary $\PP^3$ induces a K3 transition to $Y_0$, the double cover of $\PP^3$ branched in a surface $B$ of degree 10. The exceptional locus of $\pi_P$ is the quartic surface $S$, and $S$ touches $B$ in a curve of degree 20. Moreover $Y_0$ has 120 ordinary double points along $S$. Moving $B$ to general position gives a smoothing of $Y_0$.
\end{exa}

\subsection{Smoothing elliptic singularities}\label{sec!smoothings}
A necessary condition for existence of a given K3 transition is that the elliptic singularity is locally smoothable. Let $\Ce_S\subset\AA^{n+1}$ be the affine cone over a quasismooth K3 surface $S\subset\PP^n(\al)$.
\begin{conj}\label{conj!smoothings} The cone $\Ce_S$ is smoothable if and only if $S$ is the elephant section of a quasismooth Fano 3-fold $W$.
\end{conj}
One direction is easy. Suppose $S$ is an anticanonical elephant section of the Fano 3-fold $W$, and let $\Ce_W\subset\AA^{n+2}$ be the affine cone over $W$. Then by ``sweeping out the cone'' (cf.~\cite{Pinkham}), we obtain a smoothing of $\Ce_S$.
As evidence for this conjecture we present the following theorem for general smooth K3 surfaces, which is proved in \cite{Coughlan-Sano}:
\begin{thm} Let $S\subset\PP^g$ be a general smooth K3 surface of genus $g$. The affine cone $\Ce_S$ over $S$ is smoothable if and only if $g\le10$ or $g=12$.
\end{thm}
This Theorem should be compared with Pinkham's study of elliptic surface singularities \cite{Pinkham}, which are smoothable for degrees $\le9$. In applications to surface theory, the elliptic surface singularities of degree $\le3$ tend to be most useful \cite{Persson}. Similarly, here we only need K3 surfaces whose affine cone is known to be smoothable. Moreover, we need to consider the existence of global smoothings of the elliptic singularities, not just local smoothings. In this article, we mostly circumvent the issue of local versus global, because global smoothings exist more or less automatically for 3-folds of codimension $\le3$ (see Section \ref{sec!smoothings-explained}). Thus the absence of a proof of Conjecture \ref{conj!smoothings} does not pose any serious obstacles to our entertainment for now.

In Section \ref{sec!Tom-Jerry} we give two examples of quasismooth simple elliptic singularities in codimension 4, both of which have at least two topologically distinct smoothing components. Both examples are affine cones over a K3 surface that is birational to the complete intersection of a quadric and a cubic in $\PP^4$, on which we impose certain special configurations of curves. We expect many more instances of this type of behaviour, which will be investigated more fully in \cite{Brown-Coughlan}.
\begin{thm}\label{thm!obstructed} There are certain quasismooth K3 surfaces $S$ embedded in codimension 4, for which the projective cone over $S$ deforms in two different ways to topologically distinct quasismooth Fano 3-folds. These deformations induce two smoothings of the affine cone over $S$ with topologically distinct Milnor fibres.
\end{thm}
Theorem \ref{thm!obstructed} connects different Tom and Jerry families in the Hilbert scheme of Fano 3-folds \cite{BKR}.

\subsection{The web of canonical 3-folds}
We consider the directed graph whose nodes are families of canonical 3-folds, with an arrow joining two nodes if there is a K3 transition between them. We call this graph the \emph{web of canonical 3-folds} under K3 transitions. Motivated by the question of whether a triple $\chi, P_2, \Be$ is realised as a canonical 3-fold, it is natural to study fundamental properties of this web, such as whether it is connected. As a prototype, we present the following theorem.
\begin{thm}\label{thm!codim-3} Of the 137 families of canonical 3-folds in codimension $\le3$, 121 form a connected graph under K3 transitions.
\end{thm}
Based on this evidence, it seems too optimistic to expect the graph to be connected using only K3 transitions. We can improve the situation somewhat by making a compromise; we often predict numerical candidate K3 transitions between two nodes even when we can not prove existence of the transition itself. Allowing five such, we obtain a connected graph with 136 nodes. Just one family is missing.

Alternatively, we can allow more general types of transition. The behaviour of such transitions is worse than for K3 transitions, but the advantage is that we know how to construct them. Allowing five transitions which contract singular rational surfaces to non-Gorenstein non-isolated singularities, we again obtain a connected graph with 136 nodes. The missing family is different from the missing family in the numerical setting.

Thus allowing numerical K3 transitions or more general transitions, we improve Theorem \ref{thm!codim-3}:
\begin{thm-prime} Allowing five additional transitions, 136 of the 137 canonical 3-folds form a connected graph under K3 transitions.
\end{thm-prime}
We make more precise remarks about this in Sections \ref{sec!missing} and \ref{sec!non-K3}. We speculate that there could be further links passing via canonical 3-folds in codimension $\ge4$, for which the Graded Ring Database \cite{GRDB} is not yet complete.

Almost all of the K3 transitions in our web are realised as type I unprojections, and we need three unprojections of type $\two_1$. The five candidate K3 transitions mentioned above would be more difficult unprojections, so-called higher type $\two$. 

As mentioned above, an essential part of the proof of Theorem \ref{thm!codim-3} is the fact that Gorenstein rings in codimension $\le3$ are well understood, so that existence of smoothings can be readily checked.
In future work \cite{Brown-Coughlan}, we will use K3 transitions to construct and study new canonical 3-folds in codimension $>3$. 

In Section \ref{sec!props} we list some useful properties of K3 transitions without proofs. These will be taken up in \cite{Brown-Coughlan}. Section \ref{sec!web} describes the web of canonical 3-folds in codimension $\le3$ and contains some easy examples and pictures proving Theorem \ref{thm!codim-3}. In Section \ref{sec!examples} we explain some of the more complicated K3 transitions required in the proof of \ref{thm!codim-3}. Section \ref{sec!Tom-Jerry} is a detailed construction proving Theorem \ref{thm!obstructed} concerning elliptic singularities with obstructed smoothings.

\subsection{Acknowledgements}
The idea of K3 transitions was inspired by a reading of \cite{C-3f}. I thank Paul Hacking, Miles Reid and Taro Sano for helpful comments. I also thank Gavin Brown and Al Kasprzyk for their comments and for providing computer code to interact with the Graded Ring Database. This work was support by the DFG through grant Hu 337-6/2.

\section{First properties of K3 transitions}\label{sec!props}
We begin by formalising the notion of a Gorenstein elliptic 3-fold singularity.
\begin{dfn} Suppose that $P\in X$ is a Gorenstein singular point and let $f\colon \Xtilde\to X$ be a resolution of $P$. Then $P\in X$ is called an \emph{elliptic} singularity if $R^2f_*\Oh_{\Xtilde}=\CC_P$, the skyscraper sheaf supported at $P$.
\end{dfn}
Next we collect together some results about the invariants of two canonical 3-folds $X$ and $Y$ that are related by a K3 transition $X\to Y$. Proofs will appear in \cite{Brown-Coughlan}.
\begin{thm}\label{thm!invts} Let $X\to Y$ be a K3 transition. Then
\begin{enumerate}
\item $K_X^3=K_Y^3+A^2\text{ where }A=\Oh_S(1)$;
\item $\chi(\omega(X))=\chi(\omega(Y))+1$;
\item $\Be(Y)=\Be(X)\cup\Be(S)$;
\item $K_X\cdot c_2(X)=K_Y\cdot c_2(Y)+\delta(S)$, where $\delta(S)=24-\sum_{\Be(S)}\frac{r^2-1}{r}$.
\end{enumerate}
\end{thm}
A cyclic quotient singularity of type $\frac1r(a,-a)$ on $S$ induces one of type $\frac1r(a,-a,-1)$ in the basket of $Y_0$ and $Y$. This explains the union of baskets in part 3 of the theorem.
\begin{exa}\label{exa!basket} There is a K3 transition from $X_7\subset\PP(1,1,1,1,2)$ to $Y_8\subset\PP(1,1,1,2,2)$ via the K3 surface $S_6\subset\PP(1,1,2,2)$. The baskets of $X$ and $Y$ are $\Be(X)=\{\frac12(1,1,1)\}$ and $\Be(Y)=\{4\times\frac12(1,1,1)\}$. The three extra $\frac12(1,1,1)$ points on $Y$ are induced by the three $A_1$ singularities on $S$.
\end{exa}
For a polarised variety $V,A$, we write $P_V(t)=\sum_{n\ge0}h^0(V,nA)t^n$ for the Hilbert series of $V$. In terms of Hilbert series, the theorem says
\begin{cor}\label{cor!hilbert} Suppose $X$ and $Y$ are as above. Then
\[P_X(t)=P_Y(t)+\frac{t}{1-t}P_S(t),\]
where $P_S$ is the Hilbert series of the K3 surface $S$ polarised by $A$.
\end{cor}
The corollary can be proved by direct comparison of terms in orbifold Riemann--Roch for canonical 3-folds and quasismooth K3 surfaces \cite{YPG}. A more general version appears in \cite[\S2.7]{Papadakis-Reid} independent of Theorem \ref{thm!invts}. The Corollary is used throughout our proof of Theorem \ref{thm!codim-3}.

Following on from Section \ref{sec!smoothings}, suppose that $S$ is the anticanonical surface section of a Fano 3-fold $W$, smoothed by sweeping out the cone.
\begin{thm}\label{thm!euler} Let $X\to Y$ be a K3 transition in $S$. Then
\[e(X)-e(Y)=e(W)-2e(S)+2N,\]
where $N$ is the number of ordinary double points on $Y_0$.
\end{thm}
This uses the fact that the Milnor fibre is homotopic to $W-S$, combined with the standard computation of the change in Euler number coming from the ordinary double point conifold transitions.

\section{Investigating the web}\label{sec!web}
In this section, we prove Theorem \ref{thm!codim-3}. The full web is documented in \cite{tables}, and requires approximately 130 computations proving existence of each of the K3 transitions.
\subsection{Algorithm}
We start from the list of all codimension $\le3$ canonical 3-folds.
For each pair $X$, $Y$ with $p_g(X)=p_g(Y)+1$, use Corollary \ref{cor!hilbert} to compute the Hilbert series of the candidate K3 giving a possible transition $X\to Y$. 
Next search the Graded Ring Database \cite{GRDB} for existence of a K3 surface matching that Hilbert series. A posteriori, we only find K3 surfaces of codimension $\le3$. This raw data can be found in \cite{tables}.

\begin{exa} The link constructed in Example \ref{ex!intro} is a K3 transition from $X_6\subset\PP^4$ to $Y_{10}\subset\PP(1,1,1,1,5)$ via $S_4\subset\PP^3$. By Corollary \ref{cor!hilbert}, we have $P_X(t)=P_Y(t)+\frac t{1-t}P_S(t)$, where
\[P_X(t)=\frac{1-t^6}{(1-t)^5},\ P_Y(t)=\frac{1-t^{10}}{(1-t)^4(1-t^5)},\text{ and }P_S(t)=\frac{1-t^4}{(1-t)^4}.\]
\end{exa}

The next step is to prove existence or nonexistence of each numerical K3 transition. We use unprojection to contract the K3 surface, and Gorenstein formats to show that a smoothing exists. We give full specifications for each K3 transition appearing in our web in \cite{tables}.

\subsection{Failing K3 transitions}
We give two examples of K3 transitions which fail to exist, even though Corollary \ref{cor!hilbert} suggests that they might. This gives some simple checks to remove bogus links from the raw data.
\begin{exa} $121\to 1$. The target of the K3 transition would be $Y_6\subset\PP^4$ containing a K3 surface $S_{2,2,2}\subset\PP^5$. This is not possible because the embedding dimension of $S$ is too high.
\end{exa}
\begin{exa} $35\to41$. The target would be $Y_{3,10}\subset\PP(1,1,1,2,2,5)$ containing $S_6\subset\PP(1,1,2,2)$. Thus $S$ is a codimension 3 complete intersection with equation degrees $(1,5,6)$ in $\PP(1,1,1,2,2,5)$. There is no irreducible $Y$ containing $S$, because all cubics in the ideal of $S$ are divisible by the ideal generator of degree 1.
\end{exa}
In \cite{tables}, we describe five numerical K3 transitions which are higher type $\two$ unprojections. We do not yet have a construction for them, but they not fail for any obvious reason.

\subsection{Drawing the web}
We draw a directed graph, where each node is a canonical 3-fold, labelled by its ID from the Graded Ring Database. An arrow connects two nodes if there is a K3 transition from the tail to the head. Thus Examples \ref{ex!intro} and \ref{exa!basket} appear as $1\to6$ and $2\to3$ in Figure \ref{fig!hypersurface-web}. For simplicity, we draw only snapshots of the full web, highlighting the main features. We draw sufficiently many transitions to obtain a spanning tree for each connected component of the graph. There are many more K3 transitions than those appearing in the Figures, but these are redundant and so we eliminate them from our enquiries.

\subsubsection{Hypersurfaces}
\begin{figure}[ht]
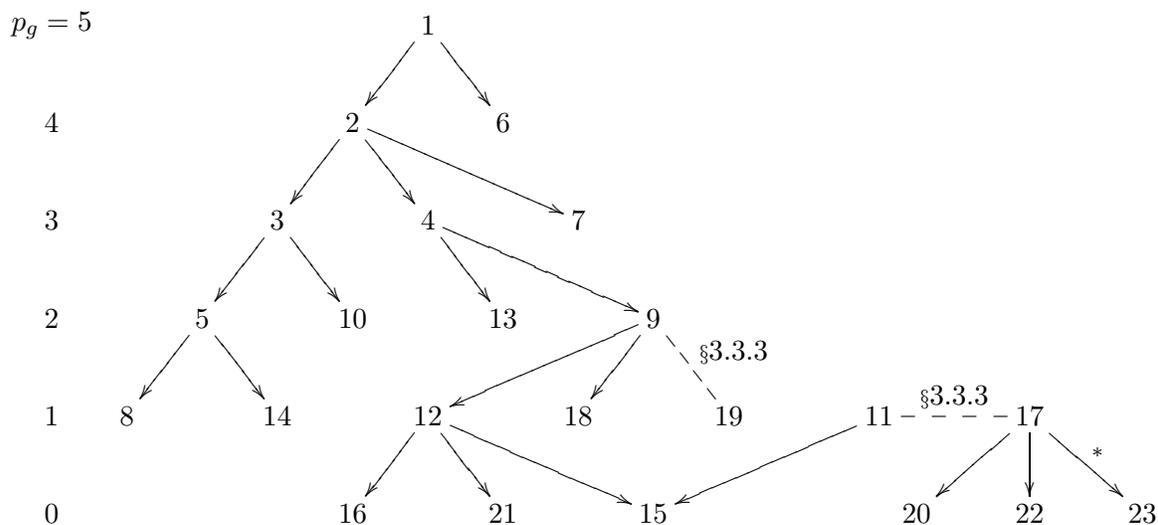

\[\xygraph{
!{<0cm,0cm>;<1cm,0cm>:<0cm,1.3cm>::}
!{(-5,5)}*+{{p_g=5}}
!{(0,5)}*+{{1}}="1"
!{(-5,4)}*+{{4}}
!{(-1,4)}*+{{2}}="2"
!{(1,4)}*+{{6}}="6"
!{(-5,3)}*+{{3}}
!{(-2,3)}*+{{3}}="3"
!{(0,3)}*+{{4}}="4"
!{(2,3)}*+{{7}}="7"
!{(-5,2)}*+{{2}}
!{(-3,2)}*+{{5}}="5"
!{(-1,2)}*+{{10}}="10"
!{(1,2)}*+{{13}}="13"
!{(3,2)}*+{{9}}="9"
!{(-5,1)}*+{{1}}
!{(-4,1)}*+{{8}}="8"
!{(-2,1)}*+{{14}}="14"
!{(0,1)}*+{{12}}="12"
!{(2,1)}*+{{18}}="18"
!{(4,1)}*+{{19}}="19"
!{(6,1)}*+{{11}}="11"
!{(8,1)}*+{{17}}="17"
!{(-5,0)}*+{{0}}
!{(-1,0)}*+{{16}}="16"
!{(1,0)}*+{{21}}="21"
!{(3,0)}*+{{15}}="15"
!{(6.5,0)}*+{{20}}="20"
!{(8,0)}*+{{22}}="22"
!{(9.5,0)}*+{{23}}="23"
"1":"2"
"1":"6"
"2":"3"
"2":"4"
"2":"7"
"3":"10"
"3":"5"
"4":"9"
"4":"13"
"5":"8"
"5":"14"
"9":"12"
"9":"18"
"11":"15"
"11"-@{--}"17"^{\S\ref{sec!missing}}
"9"-@{--}"19"^{\S\ref{sec!missing}}
"12":"16"
"12":"21"
"12":"15"
"17":"20"
"17":"22"
"17":"23"^{*}
}\]\caption{Canonical 3-fold hypersurface web}\label{fig!hypersurface-web}
\end{figure}

In Figure \ref{fig!hypersurface-web}, we draw the subgraph consisting of all 23 nodes corresponding to hypersurfaces. Nodes with equal latitude have equal $p_g$, ranging from $p_g=5$ for node $1$ at the top, down to $p_g=0$ at the bottom. The K3 transitions become more sparse, and the transitions more difficult, towards the bottom right of the hypersurface web. Thus, these are the 3-fold hypersurfaces which are the most special. For example, node 23 is Fletcher's famous example $X_{46}\subset\PP(4,5,6,7,23)$ with $p_g=0$ and $K^3=\frac1{420}$. All the information necessary to produce Figure \ref{fig!hypersurface-web} is in \cite{tables}, which also gives the K3 surface corresponding to each transition.

The hypersurface web has three connected components. The dashed links $9{--}19$ and $11{--}17$ are not K3 transitions; see Section \ref{sec!missing}.
Transition $17\xrightarrow{*}23$ needs a type II unprojection, and this is constructed in Section \ref{sec!starred}.

\begin{exa}\label{ex!AxminusBy} We construct K3 transition $2\to3$ from Example \ref{exa!basket}. Let $Y$ be a special hypersurface of degree 8 with equation
\[Y\colon(x_3f_7-y_2g_6=0)\subset\PP(1,1,1,2,2),\]
where $x_1,x_2,x_3,y_1,y_2$ are coordinates of the ambient space, and $f$, $g$ are general homogeneous forms of respective degrees $7$, $6$. Now, $Y$ contains $S_6\subset\PP(1,1,2,2)$ embedded as the complete intersection $Y\cap(x_3=g=0)$. We see that $Y$ has 21 ordinary double points along $S$ at points where $x_3=y_2=f=g=0$, and a $\frac12(1,1,1)$ point over each of the three $A_1$ singularities of $S$.

We introduce the unprojection variable $s=\frac fg=\frac{y_2}{x_3}$ to contract $S$ on $Y$ to $P$ in the complete intersection $X$ with equations
\[X\colon (sg=f,\ sx_3=y_2)\subset\PP(1,1,1,1,2,2).\]
The second equation defining $X$ eliminates coordinate $y_2$, and so $X$ is actually a hypersurface of degree $7$ in $\PP(1,1,1,1,2)$.  The tangent cone to $X$ at $P_s$ is $g(x_1,x_2,x_3,y_1)$. This is weighted homogeneous of degree 6 if we consider $x_3$ to have degree 2, since $x_3=y_2$ near $P_s$. We can smooth the singularity at $P_s$ because $X$ is a hypersurface. All but one of the K3 transitions appearing in Figure \ref{fig!hypersurface-web} are applications of the same trick used in this example.
\end{exa}
\subsubsection{Codimension $\ge2$}
\begin{figure}[ht]
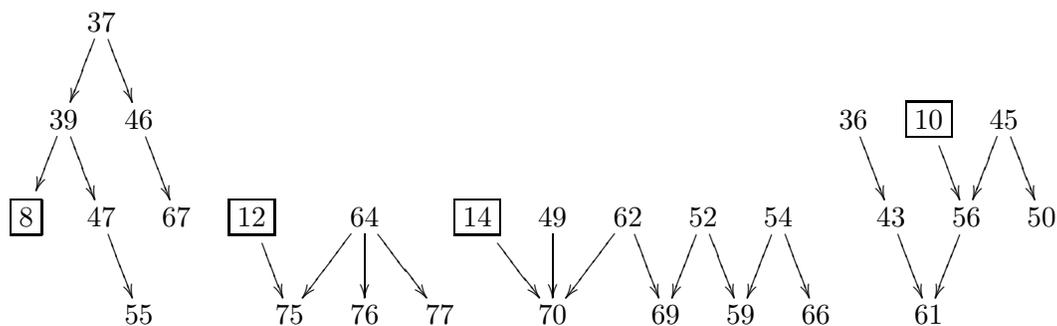

\[\xygraph{
!{<0cm,0cm>;<1cm,0cm>:<0cm,1.3cm>::}
!{(-2,3)}*+{{37}}="37"
!{(-2.5,2)}*+{{39}}="39"
!{(-1.5,2)}*+{{46}}="46"
!{(-3,1)}*+{\framebox{8}}="8"
!{(-2,1)}*+{{47}}="47"
!{(-1,1)}*+{{67}}="67"
!{(-1.5,0)}*+{{55}}="55"
!{(0,1)}*+{\framebox{12}}="12"
!{(1.5,1)}*+{{64}}="64"
!{(0.5,0)}*+{{75}}="75"
!{(1.5,0)}*+{{76}}="76"
!{(2.5,0)}*+{{77}}="77"
!{(3,1)}*+{\framebox{14}}="14"
!{(4,1)}*+{{49}}="49"
!{(5,1)}*+{{62}}="62"
!{(6,1)}*+{{52}}="52"
!{(7,1)}*+{{54}}="54"
!{(4,0)}*+{{70}}="70"
!{(5.5,0)}*+{{69}}="69"
!{(6.5,0)}*+{{59}}="59"
!{(7.5,0)}*+{{66}}="66"
!{(8,2)}*+{{36}}="36"
!{(9,2)}*+{\framebox{10}}="10"
!{(10,2)}*+{{45}}="45"
!{(8.5,1)}*+{{43}}="43"
!{(9.5,1)}*+{{56}}="56"
!{(10.5,1)}*+{{50}}="50"
!{(9,0)}*+{{61}}="61"
"37":"39"
"37":"46"
"39":"8"
"39":"47"
"46":"67"
"47":"55"
"12":"75"
"64":"75"
"64":"76"
"64":"77"
"14":"70"
"49":"70"
"62":"70"
"62":"69"
"52":"69"
"52":"59"
"54":"59"
"54":"66"
"36":"43"
"43":"61"
"56":"61"
"10":"56"
"45":"56"
"45":"50"
}\]\caption{Snapshots of the canonical 3-fold web in codimension $\le3$}\label{fig!codim-2-web}
\end{figure}

We proceed inductively on codimension. Of the 59 canonical 3-folds in codimension $2$, 24 have at least one K3 transition linking them directly to a hypersurface. A further 11 only have K3 transitions to other 3-folds of codimension 2, but after at most two such, we get a link to a hypersurface, e.g.~complete intersection $46$ is linked to hypersurface $8$ via two transitions to 3-folds of codimension 2: $46\gets37\to39\to8$. For the remaining cases, the link to a hypersurface requires several steps, and we refer to  Figures \ref{fig!codim-2-web} and \ref{fig!extra-links} for details. In each figure, hypersurface nodes are boxed.

Of the 55 canonical 3-folds in codimension 3, all but four have a K3 transition direct to a codimension 2 complete intersection, or after a single transition via another canonical 3-fold of codimension 3.

\subsubsection{Missing links}\label{sec!missing}
There are 16 nodes which are not connected to the main connected component of the web of canonical 3-folds. Of these, 12 lie in two smaller connected graphs with seven and five nodes respectively. These two smaller components appear in Figure \ref{fig!extra-links} below.
The remaining four disconnected nodes are in codimension 3; they are 113, 116, 136 and 137.

\begin{figure}[ht]
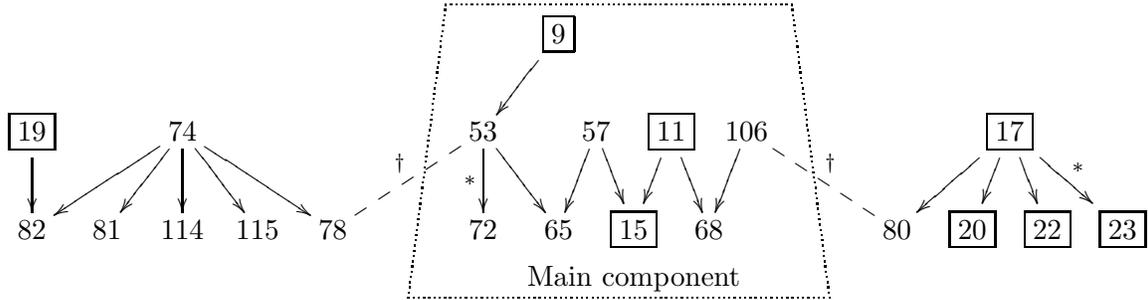

\[
\xygraph{
!{<0cm,0cm>;<1cm,0cm>:<0cm,1.3cm>::}
!{(-1,2)}*+{\framebox{9}}="9"
!{(-8,1)}*+{\framebox{19}}="19"
!{(-6,1)}*+{{74}}="74"
!{(-2,1)}*+{{53}}="53"
!{(-0.5,1)}*+{{57}}="57"
!{(-8,0)}*+{{82}}="82"
!{(-7,0)}*+{{81}}="81"
!{(-6,0)}*+{{114}}="114"
!{(-5,0)}*+{{115}}="115"
!{(-4,0)}*+{{78}}="78"
!{(-2,0)}*+{{72}}="72"
!{(-1,0)}*+{{65}}="65"
!{(0,0)}*+{\framebox{15}}="15"
!{(0,-0.5)}*+{\txt{Main component}}
!{(5,1)}*+{\framebox{17}}="17"
!{(1.5,1)}*+{{106}}="106"
!{(0.5,1)}*+{\framebox{11}}="11"
!{(1,0)}*+{{68}}="68"
!{(3.5,0)}*+{{80}}="80"
!{(4.5,0)}*+{\framebox{20}}="20"
!{(5.5,0)}*+{\framebox{22}}="22"
!{(6.5,0)}*+{\framebox{23}}="23"
!{(-3,-0.7)}="left-bottom"
!{(-2.5,2.3)}="left-top"
!{(2.6,-0.7)}="right-bottom"
!{(2.1,2.3)}="right-top"
"left-top"-@{..}"left-bottom"
"right-top"-@{..}"right-bottom"
"left-top"-@{..}"right-top"
"left-bottom"-@{..}"right-bottom"
"19":"82"
"74":"82"
"74":"81"
"74":"78"
"74":"114"
"74":"115"
"53"-@{--}"78"_{\dag}
"53":"72"_{*}
"53":"65"
"9":"53"
"57":"65"
"57":"15"
"17":"20"
"17":"22"
"17":"23"^{{*}}
"17":"80"
"106"-@{--}"80"^{{\dag}}
"106":"68"
"11":"68"
"11":"15"
}
\]\caption{Extra connected components and how they interact with the main component}\label{fig!extra-links}
\end{figure}

We have five numerical candidate K3 transitions which would connect 15 of the 16 nodes to the main connected component of the web of canonical 3-folds. These numerical transitions are
\[53\to78,\ 106\to80,\ 64\to116,\ 105\to136,\ 46\to137\]
and they are marked with a $\dagger$ in \cite{tables} and in Figure \ref{fig!extra-links}.
We do not know of any K3 transitions involving node 113, numerical or otherwise. These numerical transitions might be realised by a higher type $\two$ unprojection, but in each case, it is not easy to categorically say whether the unprojection is constructible or not, because it is a higher type $\two$ unprojection.

\subsubsection{Non-K3 transitions}\label{sec!non-K3}

As mentioned in the introduction, we can also link 15 of the 16 missing nodes to the main component by admitting as transitions, contractions of certain rational surfaces to non-Gorenstein non-isolated singularities. The following five transitions connect 15 nodes to the main component
\[12\to17,\ 12\to19,\ 64\to137,\ 77\to 113,\ 109\to136.\]
See \cite{tables} for details. Moreover, we can construct all of these, see Example \ref{ex!non-K3} for one such. The missing family is 116, which does have a non-K3 transition linking it to the rest of the web, but this would be realised by a higher type $\two$ unprojection.

\subsection{Existence of global smoothings in small codimension}\label{sec!smoothings-explained}
In the lists of Fletcher--Reid and Brown--Kasprzyk--Zhou, it is implicitly proven that the general 3-fold in any given family is quasismooth. Thus, once we construct the unprojection $Y_0\dashrightarrow X_0$, the smoothing of the elliptic singularity on $X_0$ is ensured automatically by the existence of formats for Gorenstein varieties in codimension $\le3$. Section \ref{sec!examples} works through several cases of this. In Section \ref{sec!Tom-Jerry}, we exhibit two simple elliptic 3-fold singularities embedded in codimension 4, both of which have an obstructed smoothing space.

\section{Some K3 transitions}\label{sec!examples}
In this section we compute some of the more complicated transitions appearing in the web. We have tried to cover most types of behaviour, and these serve as textbook examples of how the unprojections work in general.

\subsection{Pfaffian inside complete intersection}
Suppose $S\subset Y\subset\PP^5(\al)$, where $S$ is a K3 surface of codimension 3 defined by Pfaffians of the $5\times5$ skew matrix $M=(m_{ij})$ and $Y$ is a codimension 2 complete intersection in weighted projective space $\PP^5(\al)$. We write the two equations of $Y$ as linear combinations of the Pfaffians of $M$:
\[Y\colon\left(\sum_{i=1}^5a_i\Pf_i M=\sum_{i=1}^5b_i\Pf_i M=0\right)\subset\PP^5(\al),\]
with coefficients $a_i$, $b_i$ for $i=1,\dots,5$.
Then unprojecting $S$ in $Y$ gives $X$ in $\PP^6(1,\al)$ defined by the $6\times6$ Pfaffians of
\[\begin{pmatrix}
s & b_1 & -b_2 & b_3 & -b_4 & b_5 \\
  & a_1 & -a_2 & a_3 & -a_4 & a_5 \\
  & & m_{12} & m_{13} & m_{14} & m_{15} \\
  & & & m_{23} & m_{24} & m_{25} \\
  & & & & m_{34} & m_{35} \\
  & & & & & m_{45}
\end{pmatrix}.
\]
Here $s$ is the unprojection variable, and the tangent cone to $X$ at $P_s$ is given by the Pfaffians of $M$. Thus $X$ is singular at $P_s$ and the singularity is locally the cone over $S$, as anticipated. This unprojection appeared previously in \cite{Kapustki}.

In practice, some of the Pfaffians defining $X$ are redundant for degree reasons. Indeed, for K3 transition $128\to39$, we start from $S_{4,4,5,5,6}\subset Y_{6,7}\subset\PP(1,1,2,2,3,3)$, so that the degrees of the various matrix entries are
\[\deg\begin{pmatrix}
a_1 & a_2 & a_3 & a_4 & a_5 \\
b_1 & b_2 & b_3 & b_4 & b_5
\end{pmatrix}=
\begin{pmatrix}
0 & 1 & 1 & 2 & 2 \\
1 & 2 & 2 & 3 & 3
\end{pmatrix},\
\deg M=\begin{pmatrix}
1 & 1 & 2 & 2 \\
 & 2 & 3 & 3 \\
 & & 3 & 3 \\
 & & & 4
\end{pmatrix}.
\]
We choose $a_1=1$ because it has degree zero. Then row operations allow us to fix $b_1=0$. This introduces redundancy to the seven Pfaffians of the big matrix, eliminating two of them. The remaining five equations fit together as the $4\times4$ Pfaffians of the following $5\times5$ skew matrix
\[N=\bigwedge^2
\begin{pmatrix}
0 & m_{12} & m_{13} & m_{14} & m_{15} \\
s & -a_2 & a_3 & -a_4 & a_5
\end{pmatrix}
+
\begin{pmatrix}
b_2 & b_3 & b_4 & b_5 \\
  & m_{23} & m_{24} & m_{25} \\
  & & m_{34} & m_{35} \\
  & & & m_{45}
\end{pmatrix},\]
defining $X\subset\PP(1,1,1,2,2,3,3)$, where the first summand denotes the skew $5\times5$ matrix whose $ij$th entry is the signed $ij$-minor of the displayed $2\times4$-matrix.

We construct a global smoothing of $P$ in $X$. First choose a general quasismooth canonical 3-fold $X'$ in the same family as $X$, defined by the Pfaffians of a skew matrix $N'$ with general entries. Then take the family $\mathcal{X}\to\Delta$ defined in $\PP^6(1,\al)\times\Delta$ by the Pfaffians of $N+tN'$, where $t$ is the parameter in the small disc $\Delta$. The general fibre of $\mathcal{X}$ is quasismooth.

\subsection{The star transitions}\label{sec!starred}
We explain the three type $\two_1$ unprojections marked with a $*$.

\subsubsection{Transition $17\to23$: $S_{20}\subset\PP(2,5,6,7)$ inside $Y_{46}\subset\PP(4,5,6,7,23)$}
This is a classic type $\two_1$ unprojection, related to the elliptic involutions of \cite[\S\S4,7]{CPR}, to which we refer for more details. After coordinate changes, the general map
\[\Phi\colon\PP(2,5,6,7)\to\PP(4,5,6,7,23)\]
is given by $(\al,b,c,d)\mapsto(a=\al^2,b,c,d,e=\al g)$ where $g(\al,b,c,d)$ is a weighted homogeneous polynomial of degree 21 involving only even powers of $\al$. Let $S_{20}\subset\PP(2,5,6,7)$ be defined by equation $p+q\al$, where $p$, $q$ are of degree $20$, $18$ and again, both involve only even powers of $\al$. Write $\Phi(S)=\bar{S}\subset Y$. The free resolution of $\Oh_{\bar{S}}$ as an $\Oh_\PP$-module is
\[\Oh_{\bar{S}}\xleftarrow{(1,\al)}\Oh\oplus\Oh(-2)\xleftarrow{M} \Oh(-20)\oplus\Oh(-22)\oplus\Oh(-23)\oplus\Oh(-25)\xleftarrow{M'}
\Oh(-43)\oplus\Oh(-45)\leftarrow0\]
where
\[M=\begin{pmatrix}
p & e & aq  & ag \\
q & g & p & e
\end{pmatrix}\text{ and }
M'=\begin{pmatrix}
e & ag \\
-p & -q \\
g & e \\
-q & -p
\end{pmatrix}.\]
Now $\bar{S}$ is defined by the $2\times2$ minors of $M$, and in fact, only $M_{11}$, $M_{13}$, $M_{14}$ and $M_{24}$ are needed to generate the ideal of all minors. Thus $Y$ is defined by a weighted linear combination of these four minors, and for degree reasons, the respective coefficients must be $b,c,0,\lambda$, where $\lambda$ is a constant. The equations of the unprojected 3-fold $X$ fit into $4\times4$ Pfaffians of the skew matrix
\[\begin{pmatrix}
q & g & p & e \\
  & s_0 & \lambda & -s_1 \\
  & & s_1 & c \\
  & & & b+as_0
\end{pmatrix},\]
where $s_0$, $s_1$ are the adjoined unprojection variables of degrees $1$ and $3$ respectively.
We see that $\Pf_1$, $\Pf_3$ eliminates $c$, $e$ respectively and the other 3 Pfaffians are equivalent after this elimination process. Thus the equation defining $X_{21}\subset\PP(1,3,4,5,7)$ is
\[ps_0+qs_1=\lambda g(a,b,c,d),\]
where $\lambda c=s_0(as_0+b)-s_1^2$, and $\lambda e=q(as_0+b)+ps_1$. The tangent cone to $X$ at $P_{s_0}$ is $p+qs_1/s_0$, which is a weighted homogeneous hypersurface singularity of degree $18$ and weights $(2,3,4,9)$. A smoothing of $P_{s_0}$ on $X$ is obtained simply by varying the hypersurface to a general one.
\subsubsection{Transition $73\to21$: $S_{18}\subset\PP(2,3,4,9)$ inside
$Y_{30}\subset\PP(2,3,4,5,15)$}
In this example we describe a toric construction of a K3 transition from the top down. We can also describe it from the bottom up as in $17\to23$ above. We refer to \cite{HA} for notation and conventions concerning toric varieties via homogeneous coordinate rings and weighted blowups. Let $\FF$ be the toric variety with weight matrix
\[\begin{matrix}
& a & b & c & d & e & f & u \\
\hline
\phantom{-}H & 1 & 2 & 3 & 4 & 5 & 10 & 0 \\
-E & 0 & 2 & 3 & 4 & 6 & 9 & -1
\end{matrix}\]
that is, the matrix defines an action of $(\CC^*)^2$ on $\CC^7$, and $\FF$ is the 5-fold GIT quotient of $\CC^7$ by this action. The irrelevant ideal is $(a,b,c,d,e,f)\cap(a,u)\cap(b,c,d,e,f)$. According to \cite{HA}, $\FF$ is the weighted blowup of $\PP(1,2,3,4,5,10)$ at $P_a$ with weights $(2,3,4,6,9)$. Strictly speaking, we should write the weight matrix in its reduced form, but we do not do this here. Moreover, $\Pic\FF$ is generated by $H$, $E$ corresponding to the two factors of the $(\CC^*)^2$-action.

Let $\Xtilde$ be a complete intersection in $\FF$ of a general element from each of the linear systems $|6H-6E|$ and $|20H-18E|$. Thus, $\Xtilde$ is defined by two equations of the form $ae=F$, $f^2=G$ for appropriate polynomials $F$, $G$, and $\Xtilde$ is quasismooth in $\FF$. Define $X=\Proj\bigoplus_{n\ge0}H^0(\Xtilde,nH)$.

\begin{lemma} The image of $\fie\colon\Xtilde\to X$ is $X_{6,20}\subset\PP(1,2,3,4,5,10)$ and $X$ has a singular point at $P_a$, which is locally the cone over the K3 surface $S_{18}\subset\PP(2,3,4,9)$.
\end{lemma}
\begin{pf} Setting $u=0$ computes $\Xtilde\cap E$, and from the non-reduced weight matrix, we see that this intersection is a K3 surface $S_{6,18}\subset\PP(2,3,4,6,9)$ in $\Xtilde$ with normal bundle $N_{S/{\Xtilde}}=\Oh_S(-E)=\Oh_S(-1)$. Moreover, the coordinate of degree $6$ is eliminated near $P_a$ using the equation $ae=F$.
We already know that $\FF$ is the blowup of $P_a$ in $\PP(1,2,3,4,5,10)$ and thus $\fie$ contracts $S$ to $P_a$ on $X$.$\qed$
\end{pf}

By the adjunction formula, $K_{\Xtilde}=(H-E)_{\Xtilde}$. Now, taking $\Proj\bigoplus_{n\ge0}H^0(\FF,n(H-E))$ induces a map $\FF\to\PP(2,3,4,5,6,15)$ given in coordinates by $(b,c,d,x=ue,y=ae,z=ef)$. Thus the canonical model of $\Xtilde$ is the hypersurface $Y_{30}\subset\PP(2,3,4,5,15)$ with equation
\[z^2=[e^2G],\]
where the displayed equation just expresses $f^2=G$ in terms of the new coordinates. Again, the coordinate $y=ae$ of weight 6 is eliminated via $ae=F$, the other defining equation for $\Xtilde$.
\subsubsection{Transition $53\to72$: $S_{12}\subset\PP(2,2,3,5)$ inside $Y_{12,14}\subset\PP(2,3,4,4,5,7)$}
This transition can also be computed as a type $\two$ unprojection, or via the complete intersection $\Xtilde$ of one general element from each of the linear systems $|6H-4E|$, $|12H-12E|$ in the toric 5-fold $\FF$ with weight matrix
\[\begin{matrix}
& a & b & c & d & e & f & u \\
\hline
\phantom{-}H & 1 & 2 & 3 & 3 & 4 & 4 & 0 \\
-E & 0 & 2 & 2 & 3 & 4 & 5 & -1
\end{matrix}\]
We leave the details to the reader.

\subsubsection{A non-K3 transition $12\to17$: $S_{14}\subset\PP(1,3,4,7)$ inside $Y_{21}\subset\PP(1,3,4,5,7)$}\label{ex!non-K3}
The unprojection itself is an application of the standard hypersurface trick of Example \ref{ex!AxminusBy}. The main point is that $S$ is a rational surface which is not quasismooth. It has a singularity of the form $(x^2=yz)\subset\frac17(1,3,4)$. Thus $S$ is contracted to a nonisolated singularity on $X_{16}\subset\PP(1,2,3,4,5)$ centred at the point of Gorenstein index 2 on the ambient space. The smoothing of $X$ is clear, since it is a hypersurface.

\section{An elliptic singularity with two smoothing components}\label{sec!Tom-Jerry}
In this section we produce examples of quasismooth K3 surfaces embedded in codimension 4 whose affine cone has obstructed smoothing. We assume familiarity with Tom and Jerry unprojection \cite{BKR}, and we refer also to \cite{Brown-Georgiadis} for similar constructions in an unrelated context.

Consider the family of quasismooth K3 surfaces $S$ of genus 4 with basket $2\times\frac12(1,1)$, embedded in $\PP^6(1,1,1,1,1,2,2)$. The general such $S$ is birational to the complete intersection of a quadric and a cubic in $\PP^4$ containing two disjoint $(-2)$-curves.
\begin{thm}\label{thm!tom-jerry} There are two strata of special K3 surfaces in the above family, for which the projective cone over a surface in either stratum deforms in two different ways to topologically distinct quasismooth Fano 3-folds. These deformations induce two topologically distinct smoothings of the affine cone over the K3 surface.
\end{thm}

\begin{rmk}\label{rem!dP6} The degree 6 del Pezzo surface $S$ is famously a hyperplane section of two 3-folds: $\PP^1\times\PP^1\times\PP^1$ and $V_{(1,1)}\subset\PP^2\times\PP^2$. The two 3-folds induce distinct smoothings of the affine cone over $S$ \cite{Altmann}.
\end{rmk}

Our situation is more complicated than Remark \ref{rem!dP6}, because of the presence of weighted $\CC^*$-actions. We use Tom and Jerry unprojection; Tom is related to $\PP^2\times\PP^2$ and Jerry to $\PP^1\times\PP^1\times\PP^1$.

\begin{pf}
According to \cite{BKR}, there are at least three distinct families of quasismooth Fano 3-folds of genus 4 and basket $2\times\frac12(1,1,1)$, and we refer to them as Tom$_1$, Jerry$_{34}$ and Jerry$_{45}$.
We exhibit two strata of K3 surfaces, such that the projective cone $C_pS$ over any $S$ in stratum 1 (respectively 2),  lies in the intersection $\text{Tom}^c_1\cap \text{Jerry}^c_{34}$ (respectively $\text{Tom}^c_1\cap \text{Jerry}^c_{45}$), where the superscript means we take the closure in the Hilbert scheme of Fano 3-folds.

For each $S$ in a given stratum, we exhibit deformations of $C_pS$ to two quasismooth Fano 3-folds. Following \cite{BKR}, the node counts can be computed in terms of certain Chern classes, and these show that these Fano 3-folds are topologically distinct. Moreover, the deformation of $C_pS$ descends to the affine cone over $S$, and thus Theorem \ref{thm!tom-jerry} is proved.
$\qed$
\end{pf}

We sketch the construction of each stratum below, and in Section \ref{sec!Jerry34} we describe Jerry$_{34}$ in more detail, as it is the most interesting case.

\paragraph{Stratum 1: Tom$_1$ and Jerry$_{34}$} Start from the special complete intersection $S''_{2,3}$ in $\PP^4$ with coordinates $a,b,c,d,e$ defined by
\[S''\colon\begin{pmatrix} a & b & c \\ L & M & N \end{pmatrix}
\begin{pmatrix} a \\ -b \\ d\end{pmatrix}=0\subset\PP^4\]
where $L$ is a general quadric, and $M$, $N$ are quadrics in the ideal $(a,b,c)$.
Then $S''$ contains two lines $D_1$, $D_2$ with defining ideals $I_{D_1}=(a,b,c)$, $I_{D_2}=(a,b,d)$. The $D_i$ meet in an $A_1$ singularity on $S''$. The unprojection of $D_1$ and then $D_2$ resolves the $A_1$ singularity and contracts $D_1$ and $D_2$ to two new $\frac12(1,1)$ points on the quasismooth K3 surface $S\subset\PP(1^5,2^2)$.

Let $C_pS''\subset\PP^5$ denote the projective cone over $S''$ with cone variable $x$. There are two different deformations of $C_pS''$ to Fano 3-folds $W''\subset\PP^5$ containing two planes $\Pi_1$ and $\Pi_2$ such that $W''\cap(x=0)=S''$ and $\Pi_i\cap(x=0)=D_i$.
\begin{enumerate}
\item Tom$_1$ with 6 nodes along $\Pi_1$.
Deform $C_pS''$ to $W''_{\text{Tom}}$ by
\[W''_{\text{Tom}}\colon\begin{pmatrix} a' & b & c \\ L' & M' & N' \end{pmatrix}
\begin{pmatrix} a \\ -b' \\ d'\end{pmatrix}=0\subset\PP^5\]
where $a'=a+\al_1x$, $b'=b+\al_2x$, $d'=d+\al_3x$ for constants $\al_i$, and $L',M',N'$ can involve $x$ but we still require $M',N'$ in $(a,b,c)$. Then $W''_{\text{Tom}}$ contains two planes $\Pi_1$ and $\Pi_2$ defined by $I_{\Pi_1}=(a,b,c)$ and $I_{\Pi_2}=(a,b',d')$. The two planes meet in an ordinary double point $P_e$ on $W''$, and $W''$ has a further 6 ordinary double points on each plane, making a total of 13.

The unprojection of $\Pi_2$ in $W''$ performs a small resolution of $P_e$, and yields $W'$ containing the birational transform of $\Pi_1$ in Tom$_1$ format. Moreover, $W'$ has 6 nodes along $\Pi_1$, and the unprojection of $\Pi_1$ gives a quasismooth Fano 3-fold $W\subset\PP(1^6,2^2)$ with genus 4 and $2\times\frac12(1,1,1)$ points by the general theory of unprojection \cite{Papadakis-Reid}.

\item Jerry$_{34}$ with 8 nodes along $\Pi_1$. See Section \ref{sec!Jerry34}.

\end{enumerate}

\paragraph{Stratum 2: Tom$_1$ and Jerry$_{45}$} We start from
\[S''\colon\begin{pmatrix} a & b & c \\ L & M & N \end{pmatrix}
\begin{pmatrix} a \\ -e \\ d\end{pmatrix}=0\subset\PP^4,\]
where $L$, $M$, $N$ are all in the ideal $(a,b,c)$. Thus $S''$ contains two disjoint lines $D_1$ and $D_2$ with ideals $I_{D_1}=(a,b,c)$ and $I_{D_2}=(a,d,e)$. The unprojection of $D_1$ and $D_2$ gives a quasismooth K3 surface $S\subset\PP(1^5,2^2)$ as before. The two deformations of $C_pS''$ are
\begin{enumerate}
\item Tom$_1$ with 6 nodes on $\Pi_1$. Deform $C_pS''$ to
\[W''_{\text{Tom}}\colon\begin{pmatrix} a' & b & c \\ L' & M' & N' \end{pmatrix}
\begin{pmatrix} a \\ -e \\ d\end{pmatrix}=0\subset\PP^5\]
where $L'$ is not necessarily in $(a,b,c)$ but $M'$ and $N'$ are. The two planes meet in an ordinary double point $P_x$, and $I_{\Pi_1}=(a,b,c)$ and $I_{\Pi_2}=(a,d,e)$.

This is the same as Tom$_1$ in Stratum 1 above. The unprojection of $\Pi_2$ in $W''$ resolves $P_e$, and $W'$ has 6 nodes along $\Pi_1$ in Tom$_1$ format. The unprojection of $\Pi_1$ in $W'$ gives a quasismooth Fano 3-fold $W\subset\PP(1^6,2^2)$.

\item Jerry$_{45}$ with 7 nodes on $\Pi_1$. Deform $C_pS'$ to
\[W''_{\text{Jerry}}\colon\begin{pmatrix} a & b & c \\ L' & M' & N' \end{pmatrix}
\begin{pmatrix} a' \\ -e \\ d\end{pmatrix}=0\subset\PP^5\]
where $L'$, $M'$, $N'$ remain in the ideal $(a,b,c)$. The two planes are disjoint, and $I_{\Pi_1}=(a,b,c)$ and $I_{\Pi_2}=(a',d,e)$.

The unprojection of $\Pi_2$ in $W''$ gives $W'$, with 7 nodes along $\Pi_1$ in Jerry$_{45}$ format. The unprojection of $\Pi_1$ gives quasismooth $W\subset\PP(1^6,2^2)$.
\end{enumerate}

\begin{rmk}
Note that $W_\text{Tom}$ in stratum 1 and $W_\text{Tom}$ in stratum 2 are in the same deformation family of Fano 3-folds. There is no K3 surface whose projective cone lies in the intersection $\text{Jerry}^c_{34}\cap\text{Jerry}^c_{45}$ because one variable does not appear in the matrix as a pure power cf.~\cite[\S 5]{BKR}.
\end{rmk}
\subsection{The unprojection of Jerry$_{34}$}\label{sec!Jerry34}
Deform $C_pS''$ to $W''_{\text{Jer}}$ by
\[W''_{\text{Jer}}\colon\begin{pmatrix} a & b & c \\ L' & M' & N' \end{pmatrix}
\begin{pmatrix} a \\ -b \\ d'\end{pmatrix}=0\subset\PP^5\]
where $d'=d+\be x$, $\be$ is a constant and $L',M',N'$ involve $x$ but now $L',M'$ are not in $(a,b,c)$ and $N'$ is in $(a,b,c)$. Thus the two planes in $W''_{\text{Jerry}}$ are defined by $I_{\Pi_1}=(a,b,c)$ and $I_{\Pi_2}=(a,b,d')$. This time, $W''_{\text{Jer}}$ has 8 nodes along each plane, and moreover, the two planes meet in a line, along which $W''_{\text{Jer}}$ is singular.

The unprojection of $\Pi_2$ in $W''_{\text{Jer}}$ is defined by the five $4\times4$ Pfaffians of
\begin{equation}\label{eqn!Jerry34}
\begin{pmatrix}
d' & b & a & L' \\
& a & b & M' \\
& & c & N' \\
& & & f
\end{pmatrix},
\end{equation}
where $f$ is the new variable of weight 2. We see that the above matrix is in Jerry$_{34}$ format for the ideal $I_{\Pi_1}=(a,b,c,f)$ defining the birational image of $\Pi_1$. That is, the entries in row-columns 3 and 4 are all elements of $I_{\Pi_1}$. Moreover, $W'$ has 8 nodes along $\Pi_1$. The unprojection resolves the line $\Pi_1\cap\Pi_2$ on $W''$ to a rational surface scroll $\FF_2$, then contracts the negative section to a copy of $\PP(1,1,2)$ in $W'$.

Since $L'$, $M'$ are general quadrics, and $N'$ is in the ideal $(a,b,c)$, we can write Pfaffian 4 as $L'a-M'b+N'd'\equiv A'a-B'b+C'c$ for some $A'$, $B'$ general, and $C'$ in $(a,b,d')$. The unprojection of $\Pi_1$ is then defined by the $4\times4$ Pfaffians of \eqref{eqn!Jerry34}, together with the $4\times4$ Pfaffians of the skew matrix
\[\begin{pmatrix}
c & b & a & A' \\
& a & b & B' \\
& & d' & C' \\
& & & g
\end{pmatrix}
\]
and finally the so-called long equation 
\[fg=(L',M',N')\wedge(A',B',C'),\]
where $g$ is the second unprojection variable of weight 2, and the right hand side of the long equation is a large expression involving $A'$, $B'$, $C'$, $L'$, $M'$, $N'$ as explained in \cite[\S9.2]{BKR}.

Address: Leibniz Universit\"at Hannover, Welfengarten 1, Hannover, Germany\\
email: \verb|coughlan@math.uni-hannover.de|
\end{document}